\documentclass[11pt,reqno,tbtags]{amsart}
\usepackage{amssymb}
\usepackage{url}
\usepackage[square,numbers]{natbib}

\numberwithin{equation}{section}

\allowdisplaybreaks

\newtheorem{theorem}{Theorem}[section]

\theoremstyle{definition}
\newtheorem{example}[theorem]{Example}

\theoremstyle{remark}

\newenvironment{exx}[1]
{\begin{example}[\textbf{#1}]}
{\end{example}}

\newenvironment{romenumerate}[1][0pt]{
\addtolength{\leftmargini}{#1}\begin{enumerate}
 }{\end{enumerate}}

\newcounter{oldenumi}
{\setcounter{oldenumi}{\value{enumi}}
\begin{romenumerate} \setcounter{enumi}{\value{oldenumi}}}
{\end{romenumerate}}

\newcounter{thmenumerate}

\newcounter{romxenumerate}   

\newcounter{xenumerate}   

\newcommand{\refE}[1]{Example~\ref{#1}}

\begingroup
  \divide\count255 by 60
  \multiply\count255 by -60
  \advance\count255 by \time
  \ifnum \count255 < 10 \xdef\klockan{\the\count1.0\the\count255}
  \else\xdef\klockan{\the\count1.\the\count255}\fi
\endgroup

\newcommand{\sumjn}{\sum_{j=1}^n}
\newcommand{\prodin}{\prod_{i=1}^n}

\newcommand\set[1]{\ensuremath{\{#1\}}}

\newcommand\xpar[1]{(#1)}
\newcommand\bigpar[1]{\bigl(#1\bigr)}
\newcommand\Bigpar[1]{\Bigl(#1\Bigr)}

\newcommand\lrpar[1]{\left(#1\right)}

\newcommand\lrabs[1]{\left|#1\right|}
\def\rompar(#1){\textup(#1\textup)}    
\newcommand\xfrac[2]{#1/#2}

\newcommand\parfrac[2]{\lrpar{\frac{#1}{#2}}}

\def\xexp(#1){e^{#1}}

\newcommand\ntoo{\ensuremath{{n\to\infty}}}

\newcommand\punkt[1]{\if.#1\else.\spacefactor1000\fi{#1}}
\newcommand\iid{i.i.d\punkt}    
\newcommand\ie{i.e\punkt}
\newcommand\eg{e.g\punkt}
\newcommand\viz{viz\punkt}
\newcommand\cf{cf\punkt}
\newcommand{\as}{a.s\punkt}

\newcommand\ii{\mathrm{i}}

\newcommand{\tend}{\longrightarrow}
\newcommand\dto{\overset{\mathrm{d}}{\tend}}

\newcommand\eqd{\overset{\mathrm{d}}{=}}

\newcommand\bbR{\mathbb R}
\newcommand\bbC{\mathbb C}

\newcommand\bbZ{\mathbb Z}

\newcounter{CC}
\newcounter{cc}

\renewcommand\Re{\operatorname{Re}}
\renewcommand\Im{\operatorname{Im}}

\newcommand\E{\operatorname{\mathbb E{}}}
\renewcommand\P{\operatorname{\mathbb P{}}}

\newcommand\Exp{\operatorname{Exp}}

\newcommand\ga{\alpha}
\newcommand\gb{\beta}
\newcommand\gd{\delta}
\newcommand\gD{\Delta}
\newcommand\gf{\varphi}
\newcommand\gam{\gamma}
\newcommand\gG{\Gamma}
\newcommand\gk{\varkappa}
\newcommand\gl{\lambda}

\newcommand\gth{\theta}
\newcommand\eps{\varepsilon}

\newcommand\cS{{\mathcal S}}
\newcommand\cT{{\mathcal T}}

\newcommand\tP{\widetilde P}
\newcommand\tW{{\widetilde W}}
\newcommand\hC{{\hat C}}

\newcommand\qw{^{-1}}
\newcommand\qww{^{-2}}
\newcommand\qq{^{1/2}}
\newcommand\qqw{^{-1/2}}

\newcommand\qqqq{^{1/4}}

\renewcommand{\=}{:=}

\newcommand\intoi{\int_0^1}
\newcommand\intoo{\int_0^\infty}
\newcommand\intoooo{\int_{-\infty}^\infty}
\newcommand\oi{[0,1]}

\newcommand\oooo{(-\infty,\infty)}

\newcommand\dd{\,\textup{d}}

\newcommand{\mgf}{moment generating function}
\newcommand{\chf}{characteristic function}

\newcommand\rv{random variable}

\newcommand{\B}{\mathrm{B}}
\newcommand{\N}{\mathrm{N}}
\newcommand{\U}{\mathrm{U}}

\newcommand{\ogt}{of Gamma type}
\newcommand{\mogt}{moments of Gamma type}
\newcommand{\mgfogt}{\mgf{} of Gamma type}
\newcommand\ise{_{\text{\sc ise}}} 
\newcommand\fise{f\ise}
\newcommand\muise{\mu\ise}
\newcommand{\stablex}[1]{$#1$-stable}
\newcommand\txx[1]{T_{(#1)}}
\newcommand\xxn{(x_1,\dots,x_n)}
\newcommand\XXn{(X_1,\dots,X_n)}
\newcommand\XXnx{X_1,\dots,X_n}
\newcommand\YYn{(Y_1,\dots,Y_n)}
\newcommand\YYnx{Y_1,\dots,Y_n}
\newcommand\ZZn{(Z_1,\dots,Z_n)}
\newcommand\ZZnx{Z_1,\dots,Z_n}
\newcommand\bsa{\bar S_\ga}
\newcommand\bsx[1]{\bar S_{#1}}
\newcommand\hpsi{\widehat\psi}
\newcommand\heta{\widehat\eta}
\newcommand\hf{\widehat f}
\newcommand\mm{^{(m)}}
\newcommand\mmx[1]{^{(#1)}}

\newcommand{\Levy}{L\'evy}
\newcommand\REM[1]{{\raggedright\texttt{[#1]}\par\marginal{XXX}}}

\newcommand\urladdrx[1]{{\urladdr{\def~{{\tiny$\sim$}}#1}}}

\begin{document}
\title[Further examples with moments of Gamma type]
{Further examples with moments of Gamma type}

\date{3 December, 2010; Revised 6 February 2013}

\author{Svante Janson}
\address{Department of Mathematics, Uppsala University, PO Box 480,
SE-751~06 Uppsala, Sweden}
\email{svante.janson@math.uu.se}
\urladdrx{http://www.math.uu.se/~svante/}

\subjclass[2000]{} 


\maketitle

This is an appendix to \cite{SJ244} containing further examples.
See \cite{SJ244} for notation and for examples and equations referred to
below by numbers.
See also the further references in \cite[Addendum]{SJ244}.

This appendix will probably be extended with more examples in the future.

\appendix
\setcounter{section}1
\section{Further examples}

\begin{exx}{Rayleigh distribution} \label{ERayleigh}
  The Rayleigh distribution $R$ is the chi distribution $\chi(2)$, with density 
$x e^{-x^2/2}$. This is a special case of Example 3.6, and we have
  \begin{equation}
\E R^s = 2^{s/2}\gG(s/2+1),
\qquad -2<\Re s<\infty.	
  \end{equation}

We have $\rho_+=\infty$, $\rho_-=-2$, $\gam=\gam'=1/2$, $\gd=1/2$,
$\gk=0$, $C_1=\pi\qq$.
\end{exx}

\begin{exx}{Maxwell distribution}
  The Maxwell distribution $M$ is the chi distribution $\chi(3)$, with density 
$(2/\pi)\qq x^2 e^{-x^2/2}$.
This is a another special case of Example 3.6, and we have
  \begin{equation}
\E M^s = \frac{2^{s/2}}{\gG(3/2)}\gG\Bigpar{\frac s2+\frac32}
= \frac{2^{s/2+1}}{\sqrt\pi}\gG\Bigpar{\frac s2+\frac32},
\qquad -3<\Re s<\infty.	
  \end{equation}

We have $\rho_+=\infty$, $\rho_-=-3$, $\gam=\gam'=1/2$, $\gd=1$,
$\gk=0$, $C_1=\sqrt2$.
\end{exx}

\begin{exx}{Type-2 Beta distribution}
  The type-2 Beta distribution \cite[Chapter 4]{MSH} 
has density
  \begin{equation}
	f(x)=\frac{\gG(\ga+\gb)}{\gG(\ga)\gG(\gb)}x^{\ga-1}(1+x)^{-(\ga+\gb)},
\qquad x>0,
  \end{equation}
for two parameters $\ga,\gb>0$.
A variable $X_{\ga,\gb}$ with this distribution has moments given by 
\begin{equation}
  \E X_{\ga,\gb}^s = \frac{\gG(\ga+s)\gG(\gb-s)}{\gG(\ga)\gG(\gb)},
\qquad -\ga <\Re s<\gb.
\end{equation}
A comparison with (3.1) shows that
$X_{\ga,\gb}\eqd \gG_\ga/\gG'_\gb$, with $\gG_\ga$ and $\gG'_\gb$ independent.
In particular, see Example 3.7, the $F$ distribution is of this type (up to
a constant factor):
$F_{m,n}\eqd(n/m) X_{m/2,n/2}$.

We have $\rho_+=\gb$, $\rho_-=-\ga$, $\gam=2$, $\gam'=0$, $\gd=\ga+\gb-1$,
$\gk=0$, $C_1=2\pi/(\gG(\ga)\gG(\gb))$.

Note that $1+X_{\ga,\gb}\eqd B_{\gb,\ga}\qw$, the inverse of a (usual) Beta
distributed variable, see Example 3.4. 
Thus $1+X_{\ga,\gb}$ also has moments of Gamma type,
with
\begin{equation}
\E  (X_{\ga,\gb}+1)^s = \E B_{\gb,\ga}^{-s}
=\frac{\gG(\ga+\gb)\gG(\gb-s)}{\gG(\gb)\gG(\ga+\gb-s)},
\qquad \Re s<\gb.
\end{equation}
\end{exx}

\begin{exx}{Cauchy distribution}\label{Ecauchy}
  The Cauchy distribution with density $1/(\pi(1+x^2))$, $-\infty<x<\infty$,
equals the $t$-distribution in Example 3.8 with $n=1$. Hence, if $X$ is a
random variable with a Cauchy distribution, then
$|X|\eqd|\cT_1|\eqd F_{1,1}\qq$
and $|X|$ has moments of Gamma type
\begin{equation}\label{cauchy}
  \E |X|^s=\frac1\pi
\gG\Bigpar{\frac 12+\frac s2}\gG\Bigpar{\frac 12-\frac s2}
=\frac1{\cos(\pi s/2)}
,
\qquad -1<\Re s<1.
\end{equation}
Cf.\ Example 3.19, where $A\eqd \frac2\pi\log|X|$.

We have $\rho_+=1$, $\rho_-=-1$, $\gam=1$, $\gam'=0$, $\gd=0$,
$\gk=0$, $C_1=2$.
\end{exx}

\begin{exx}{Beta product distribution}
Dufresne \cite{Dufresne:BP} has shown that 
if $a,b,c,d$ are real, then
there 
exists a probability
distribution
$G(a,c;a+b,c+d)$ on $(0,1)$ with moments
\begin{equation}\label{BP}
  \E X^s = C\frac{\gG(a+s)\gG(c+s)}{\gG(a+b+s)\gG(c+d+s)}
\end{equation}
(where necessarily $C=\gG(a+b)\gG(c+d)/\xpar{\gG(a)\gG(c)}$),
if and only if either
\begin{romenumerate}
\item $a>0$, $c>0$, $b+d>0$ and $\min(a+b,c+d)>\min(a,c)$, or 
\item 
\eqref{BP} degenerates to $\E X^s = C\xfrac{\gG(\ga+s)}{\gG(\ga+\gb+s)}$
with $\ga>0$ and $\gb\ge0$, so $X$ has a Beta distribution or $X\equiv1$.
(This degenerate case occurs if $b=0$, $d=0$, $a+b=c$ or $c+d=a$.)
\end{romenumerate}

We have $\rho_+=\infty$, $\rho_-=-\min\set{a,c}$, $\gam=\gam'=0$, $\gd=-b-d$,
$\gk=0$, $C_1=C$.

The case when all $a,b,c,d>0$ is just a product of two independent Beta
variables $B_{a,b}B_{c,d}$, 
see Example 3.4,
but there are also other possible parameter
values, for example $(2,5,8,-1)$ given in \cite{Dufresne:BP}.

Moreover, if we
allow complex parameters $a,b,c,d$, there is exactly one more case 
\cite{Dufresne:BP}, \viz{}
$a>0$, $c>0$, $b+d>0$ (entailing $\Im(b)=-\Im(d)$), and $\Re(a+b)=\Re(c+d)$.
In particular, we may take $a=c>0$ and $d=\overline b$ for any complex $b$
with $\Re b>0$. 
However, complex parameters are not included in the class of
distributions studied in \cite{SJ244}, see Remark 11.3.
\end{exx}

\begin{exx}{Density of ISE}\label{EISE}
The ISE (integrated superbrownian excursion) is a
random probability measure
introduced by Aldous \cite{AldousISE}.
It was shown in \cite{SJ185} that the ISE \as{} is absolutely continuous,
and thus has a (random) density $\fise(x)$, $x\in\oooo$. 

The ISE can be described as the occupation measure of the head of the Brownian
snake, see 
Le Gall \cite[Chapter IV]{LeGall} or Le Gall and Weill \cite{LeGallW}
for details; see also \cite[Section 4.1]{SJ176}.
Thus $\fise(x)$ is the local time of the head of the Brownian snake.
Moreover, $\fise(x)$ arises for example as a 
limit of the vertical profile of random trees, see
\cite{Marckert},  
\cite{BM}, \cite{SJ185}, \cite{SJ222} and
\cite{Drmota}. 

The distribution of $\fise(x)$ for a fixed $x$ is given by a rather
complicated formula, see \cite{BM} and \cite{SJ185}; 
in the case $x=0$ it simplifies and 
$\fise(0)\eqd 2\qqqq 3\qw S_{2/3}\qqw$, 
where $S_{2/3}$ is a
positive \stablex{2/3} variable with Laplace transform 
$\E e^{-t S_{2/3}}=e^{-t^{2/3}}$, see Example 3.10.
Thus $\fise(0)$ has moments of Gamma type with
\begin{equation}
 \E \fise(0)^s=2^{s/4}3^{-s} \frac{\Gamma(3s/4+1)}{\Gamma(s/2+1)},
\qquad -4/3<\Re s<\infty,
\end{equation}
see \cite{BM} and \cite{SJ185}.
We have $\rho_+=\infty$, $\rho_-=-4/3$, $\gam=\gam'=1/4$, $\gd=0$,
$\gk=-\frac34\log 2-\frac14\log 3$, $C_1=\sqrt{3/2}$.
\end{exx}

\begin{exx}{Average ISE}
The ISE in \refE{EISE} is a random probability measure $\muise$; taking the
expectation we obtain a deterministic probability measure $\E\muise$, which
is the distribution of a random variable $X$ that can be seen as a random
point given by a random ISE. (This is, for example, the limit distribution of
the label of a random node in a random tree under suitable assumptions and
normalizations.) $X$ has a symmetric distribution, 
and $|X|$ has moments of Gamma type with
\begin{equation}
  \E |X|^s
=
\frac{2^{3s/4}}{\sqrt\pi}
 \Gamma\Bigpar{\frac s2+\frac12}\Gamma\Bigpar{\frac s4+1},
\qquad -1<\Re s<\infty,
\end{equation}
see \cite{AldousISE}  and \cite{SJ185}.
We have $\rho_+=\infty$, $\rho_-=-1$, $\gam=\gam'=3/4$, $\gd=1/2$,
$\gk=-\frac14\log 2$, $C_1=\sqrt\pi$.
\end{exx}

\begin{exx}{Blocks in a Stirling permutation}
Let $k\ge2$ be a fixed integer. It is shown in \cite{SJ217} that
the number of blocks in a random $k$-Stirling permutation of order $n$
(see \cite{SJ217} for definitions) 
after suitable normalization converges in distribution as $\ntoo$ to a
random variable $\zeta$ with moments of Gamma type given by
\begin{equation}
  \E \zeta^s
=(s+1)!\frac{\Gamma(1+\frac1k)}{\Gamma(1+\frac{s+1}k)}
=\Gamma\Bigpar{1+\frac1k}\frac{\gG(s+2)}{\gG\bigpar{\frac sk+\frac{k+1}k}},
\qquad -2<\Re s<\infty. 
\end{equation}
As explained in \cite{SJ217}, this is actually a special case of (9.1).  

We have $\rho_+=\infty$, $\rho_-=-2$, $\gam=\gam'=\gd=(k-1)/k$, 
$\gk=\frac1k\log k$, $C_1=k^{(k+2)/2k}\gG((k+1)/k)$.
\end{exx}

\begin{exx}{Distances in a sphere}
Let $X_1$ and $X_2$ be two independent random points, uniformly distributed
in an $n$-dimensional ball of radius $a$, and let $D\=|X_1-X_2|$ be the
distance between them. (Here $n\ge1$.)
Note that $0\le D\le 2a$, so $D/2a\in\oi$.
\citet{Hammersley} showed that the density function of $D/2a$ is
\begin{equation}\label{dist_n}
f_n(\gl)=
\frac{2n\gG(n+1)\gl^{n-1}}{\gG\bigpar{\frac12n+\frac12}^2} 
  \int_{\gl}^1(1-z^2)^{(n-1)/2}\dd z
\end{equation}
and as a consequence, for $\Re s>-n$,
\begin{equation}
\E (D/2a)^s=
  \frac{n\gG(n+1)}{\gG\bigpar{\frac12n+\frac12}} 
\cdot
\frac{\gG\bigpar{\frac12n+\frac12s+\frac12}}
  {(n+s)\gG\bigpar{n+\frac12s+1}}  
\end{equation}
and, equivalently,
\begin{equation}\label{dist2}
  \begin{split}
\E D^s
&=
C(2a)^s\frac{\gG(s+n)\gG\bigpar{\frac12s+\frac12n+\frac12}}
  {\gG(s+n+1)\gG\bigpar{\frac12s+n+1}}
\\&
=C'a^s\frac{\gG(s+n)}
  {\gG\bigpar{\frac12s+\frac12n+1}\gG\bigpar{\frac12s+n+1}}
  \end{split}
\end{equation}
with $C=  \xfrac{n\gG(n+1)}{\gG\bigpar{\frac12n+\frac12}} $
and $C'=\pi\qq2^{-n}C$.
(\citet{Hammersley} did not specify the range of $s$, and presumably
intended only positive and perhaps integer values, but the formula follows
by \eqref{dist_n} for any $s$ with $\Re s>n$. Alternatively, the result
extends from positive $s$ by Theorem 2.1.)

$D$ thus has moments of Gamma type, with 
$\rho_+=+\infty$, $\rho_-=-n$, 
$\gam=\gam'=0$, $\gd=-(n+3)/2$, $\gk=\log(2a)$,
$C_1=2^{(n+1)/2}C$.

It follows from \eqref{dist2} that if $B_{n,1}$ and $B_{(n+1)/2,(n+1)/2}$
are independent Beta distributed variables, then
$\E (D/2a)^s=\E B_{n,1}^s B_{(n+1)/2,(n+1)/2}^{s/2}$, see Example 3.4, and
thus
\begin{equation}
  D\eqd 2a B_{n,1} B_{(n+1)/2,(n+1)/2}^{1/2}.
\end{equation}
Cf.\ Remark 1.5.
Note further that $B_{n,1}\eqd U^{1/n}$ where $U\sim\U(0,1)$ is uniform, see
Example 3.3, so we also have 
$  D\eqd 2a U^{1/n} B_{(n+1)/2,(n+1)/2}^{1/2}$.

Taking $n=2$ and $s=1$ in \eqref{dist2} we see that the average distance between
two random points in a circular disc of radius $a$ is
$\frac{128}{45\pi}a$; this is an old problem.

Further examples for a ball of diameter 1 (so $a=1/2$):
if $n=1$, then $\E D=1/3$ and $\E D^2=1/6$;
if $n=2$, then $\E D=64/45\pi$ and $\E D^2=1/4$;
if $n=3$, then $\E D=18/35$ and $\E D^2=3/10$.

For $s=2$, $\E D^2=2a^2n/(n+2)$, as is easily seen directly.
\end{exx}

\begin{exx}{Preferential attachment random graph}
\citet{PRR} have, motivated by the study of vertex degrees in a preferential
attachment random graph, studied a special case of the triangular urn in
Section 9 and obtained further results. (They obtain the random variable
$K_\ga$ below, for $\ga\in\set{\frac12,1,\frac32,2,\dots}$, as the limit in
distribution, after normalization,
of the degree of a fixed vertex in one of two slightly different random
 graphs.) 

In our notation, let $W_\ga$, for $\ga\ge1/2$, be the (limit) variable in
Section 9.1 with $a=2$, $c=d=1$, $w_0=1$ and $b=2\ga-1$.
Then (9.1) yields 
\begin{equation}
  \label{prr1}
\E W_\ga^s=\gG(\ga)\frac{\gG(s+1)}{\gG(s/2+\ga)},
\qquad \Re s>-1.
\end{equation}
(Theorem 2.1 implies that the condition $\ga\ge1/2$ also is necessary for
the existence of such a random variable: $\ga\in\set{0,-1,-2,\dots}$ is
clearly impossible, since then $\gG(\ga)=\infty$, and otherwise the function
in \eqref{prr1} har $\rho_-\le-1$ and $\rho_+=\infty$, and if $\ga<1/2$ it
is 0 at $s=-2\ga\in(\rho_-,\rho_+)$, which contradicts Theorem 2.1.)

\citet{PRR} choose a different normalisation, so we define
$K_\ga\=(\ga/2)\qq W_\ga$ and obtain
\begin{equation}
  \label{prr2}
\E K_\ga^s=\parfrac{\ga}2^{s/2}\frac{\gG(\ga)\gG(s+1)}{\gG(s/2+\ga)},
\qquad \Re s>-1.
\end{equation}
In particular, $K_\ga$ satifies the normalisation $E K_\ga^2=1$.

$K_\ga$ has $\rho_+=\infty$, $\rho_-=-1$ (except when $\ga=1/2$; then
$\rho_-=-2$),  
$\gam=\gam'=1/2$, $\gd=1-\ga$,
$\gk=\frac12\log\ga$, $C_1=2^{\ga-1/2}\gG(\ga)$.

\citet{PRR} show, among other things, that the random variable $K_\ga$ has
the density function
\begin{equation}
  \gk_\ga(s)=\gG(\ga)\sqrt{\frac2{\ga\pi}}\,
 e^{-x^2/2\ga}U\Bigpar{\ga-1,\frac12;\frac{x^2}{2\ga}},
\qquad x>0,
\end{equation}
where $U(a,b;z)$ denotes the confluent hypergeometric function of thesecond
kind; see \eg{} \cite[Chapter 13]{AS} or \cite{Lebedev} (where it is denoted
$\Psi$). 
This is a considerably simpler formula than the power series expansion given
in Theorem 9.1.
It would be interesting to know whether the density in Theorem 9.1 can be
expressed using hypergeometric functions also for other triangular urns.

Note the special case $\ga=1/2$; then \eqref{prr2} simplifies by the
duplication formula for the Gamma function to 
\begin{equation}
\E K_{1/2}^s=\gG(s/2+1),
\qquad \Re s>-2.
\end{equation}
showing that 
\begin{equation}
K_{1/2}\eqd T\qq \eqd 2\qqw R,  
\end{equation}
with $T\sim\Exp(1)$, see Example 3.2, 
and $R\sim\chi(2)$ (the Rayleigh distribution), see Examples 3.6 and
\ref{ERayleigh}. 
The density function of $K_{1/2}$ is thus
\begin{equation}
  \gk_{1/2}(x)=2xe^{-x^2},
\qquad x>0.
\end{equation}
\end{exx}

\begin{exx}{The maximum of \iid{} exponentials}\label{EmaxExp}

Let $(T_i)_{i=1}^\infty$ be \iid{} exponential random variables with
$T_i\sim\Exp(1)$, and let $M_n\=\max_{1\le i\le n}T_i$.
Then $e^{-T_i}\sim\U(0,1)$, and thus,  
since 
$e^{-M_n}\=\min_{1\le i\le n}e^{-T_i}$,
\begin{equation}
\P\bigpar{e^{-M_n}>x} =\P\bigpar{e^{-T_1}>x}^n=(1-x)^n,
\qquad 0<x<1,
\end{equation}
so $e^{-M_n}$ has the Beta distribution $\B(1,n)$.
  
Hence, by Example 3.4, 
\begin{equation}\label{emax}
\E{e^{s M_n}} =\E B_{1,n}^{-s}=\frac{\gG(n+1)\gG(1-s)}{\gG(n+1-s)},
\qquad \Re s<1.
\end{equation}
Hence $M_n$ has \mgf{} of Gamma type.
The special case $n=1$ gives $M_1=T_1\sim\Exp(1)$ treated in Example 3.16.

$M_n$ has $\rho_+=1$, $\rho_-=-\infty$,
$\gam=\gam'=0$, $\gd=-n$,
$\gk=0$, $C_1=n!$,
\cf{} Example 3.4 and Remark 2.8.

For an alternative proof of \eqref{emax}, note that if $\txx1<\dots<\txx n=M_n$
are $T_1,\dots,T_n$ arranged in increasing order, then it is a standard
observation
(\eg{} by regarding $T_1,\dots,T_n$ as the first points in independent
Poisson processes) that $\txx1$, $\txx2-\txx1, \dots\txx{n}-\txx{n-1}$
are independent exponential variables with
$\txx k-\txx{k-1}\sim\Exp(1/(n-k+1))$
(wih $\txx0\=0$), and hence,
for $\Re s<1$
\begin{equation}\label{apr}
  \begin{split}
  \E e^{sM_n} 
&= \prod_{k=1}^n \frac1{1-s/(n-k+1)}
= \prod_{j=1}^n \frac1{1-s/j}
= \prod_{j=1}^n \frac{j}{j-s}
\\&
=\gG(n+1)\frac{\gG(1-s)}{\gG(n+1-s)}.	
  \end{split}
\end{equation}
Note further that, as \ntoo,
\begin{equation}\label{febr}
  \begin{split}
  \E e^{s(M_n-\log n)} 
=n^{-s}\frac{\gG(n+1)}{\gG(n+1-s)}{\gG(1-s)}
\to {\gG(1-s)},
\qquad \Re s<1,	
  \end{split}
\end{equation}
and thus
\begin{equation}\label{morm}
  M_n-\log n\dto W,
\end{equation}
where $W$ has the Gumbel distribution
$\P(W\le x)=e^{-e^{-x}}$ which has \mgf{} $\E e^{sW}=\gG(1-s)$, $\Re s<1$,
see Example 3.19.
Note that it is also easy to prove \eqref{morm} directly, since, for
$x\in\bbR$ and 
$n$ large enough, 
\begin{equation}
  \P(M_n-\log n\le x) =\P(T_1\le \log n+x)^n
=\lrpar{1-\frac{e^{-x}}n}^n\to  e^{-e^{-x}}.
\end{equation}

The decomposition
\begin{equation}\label{terminalia}
  M_n=\sum_{k=1}^n \bigpar{\txx k-\txx{k-1}}
\eqd \sum_{k=1}^n \frac{1}{n-k+1}T_k
\eqd \sum_{j=1}^n \frac{1}{j}T_j
\end{equation}
shows also that 
\begin{equation}
  \E M_n = \sumjn \frac{1}j,
\end{equation}
the $n$:th harmonic number $H_n$.
We have $ \E W=-\gG'(1)=\gamma$,
Euler's gamma, and thus, since \eqref{febr} implies convergence of all
moments,
\begin{equation}\label{harm}
  H_n-\log n =\E(M_n-\log n) \to \E W=\gamma,
\end{equation}
a well-known result by \citet{Euler}.

Moreover, it follows from \eqref{febr},  or from \eqref{morm} and
\eqref{harm}, that
\begin{equation}
  M_n-\E M_n\dto W-\E W = W -\gamma,
\end{equation}
and thus by \eqref{terminalia}
\begin{equation}\label{wsum}
  \sum_{j=1}^\infty \frac{1}{j}(T_j-1)
\eqd W-\gamma,
\end{equation}
where the infinite sum converges in $L^2$ and thus \as{}
\cite[Lemma 4.16]{Kallenberg}.
\end{exx}

\begin{exx}{The largest values of \iid{} exponentials}\label{EmaxExp2}
Gen\-er\-al\-iz\-ing \refE{EmaxExp},
let $M_n\mm$ be the $m$:th largest of the 
$n$ \iid{} exponential random variables
$T_1,\dots,T_n\sim\Exp(1)$;
here $1\le m\le n$. 
(The special case $m=1$ gives $M_n$ treated in \refE{EmaxExp}.)

Let $U_i\=e^{-T_i}\sim\U(0,1)$. Then $e^{-M_n\mm}$ is the $m$:th smallest of the
\iid{} uniform $U_1,\dots,U_n$, and thus
$e^{-M_n\mm}$ has the Beta distribution $\B(m,n-m+1)$.
  
Hence, by Example 3.4, 
\begin{equation}\label{emaxm}
\E{e^{s M_n\mm}} =\E B_{m,n-m+1}^{-s}=\frac{\gG(n+1)\gG(m-s)}{\gG(m)\gG(n+1-s)},
\qquad \Re s<m.
\end{equation}
Thus $M_n\mm$ has \mgf{} of Gamma type. 

$M_n\mm$ has $\rho_+=m$, $\rho_-=-\infty$,
$\gam=\gam'=0$, $\gd=-(n-m+1)$,
$\gk=0$, $C_1=n!/(m-1)!$,
\cf{} Example 3.4 and Remark 2.8.

Alternatively,  \eqref{emaxm} can be obtained by the argument in \eqref{apr}.
Moreover, by the lack of memory for the exponential distribution,
$M_n-M_n\mm$ is independent of $M_n\mm$ and has the same distribution as
$M_{m-1}$; thus $M_n=M_n\mm+M_{m-1}'$, where $M'_{m-1}$ is a copy of
$M_{m-1}$ that is independent of $M_n\mm$; this yields 
$\E e^{sM_n} = \E e^{sM_n\mm} \E e^{sM_{m-1}}$, and \eqref{emaxm} follows
from \eqref{emax}.

As \ntoo,
\begin{equation}\label{febrm}
  \begin{split}
  \E e^{s(M_n\mm-\log n)} 
=n^{-s}\frac{\gG(n+1)}{\gG(n+1-s)}\cdot\frac{\gG(m-s)}{\gG(m)}
\to \frac{\gG(m-s)}{\gG(m)},
\qquad \Re s<m,	
  \end{split}
\end{equation}
and thus
\begin{equation}\label{mormm}
  M_n\mm-\log n\dto W\mm,
\end{equation}
where $W\mm$ has the \mgf{} of Gamma type 
\begin{equation}
  \E e^{s W\mm}= 
   \frac{\gG(m-s)}{\gG(m)},
\qquad \Re s<m.
\end{equation}

Comparing with Example 3.1, we see that
$  \E e^{s W\mm}= \E \gG_m^{-s}=\E e^{-s\log\gG_m}$
and thus $W\mm\eqd -\log\gG_m$, where $\gG_m$ has a Gamma
distribution $ \gG(m)$. 

$W\mm$ has $\rho_+=m$, $\rho_-=-\infty$,
$\gam=1$, $\gam'=-1$, $\gd=m-1/2$,
$\gk=0$, $C_1=\sqrt{2\pi}/(m-1)!$,
\cf{} Example 3.1 and Remark 2.8.

As in \eqref{terminalia}, there is a decomposition
\begin{equation}\label{terminaliam}
  M_n\mm=\sum_{k=1}^{n-m+1} \bigpar{\txx k-\txx{k-1}}
\eqd \sum_{k=1}^{n-m+1} \frac{1}{n-k+1}T_k
\eqd \sum_{j=m}^n \frac{1}{j}T_j,
\end{equation}
which shows  that 
\begin{equation}
  \E M_n = \sum_{j=m}^n \frac{1}j=H_n-H_{m-1}.
\end{equation}
Since \eqref{febrm} implies convergence of all
moments, this yields
\begin{equation}\label{harmm}
\E W\mm=\lim_{\ntoo} \E(M_n\mm-\log n) 
=\lim_{\ntoo} \bigpar{ H_n-H_{m-1}-\log n }=
\gamma-H_{m-1}.
\end{equation}

Moreover, \eqref{terminaliam},  \eqref{mormm} and \eqref{harmm} imply
\begin{equation}\label{wsumm}
  \sum_{j=m}^\infty \frac{1}{j}(T_j-1)
\eqd W\mm-\E W\mm
\eqd W\mm-\gamma+H_{m-1},
\end{equation}
where the infinite sum converges in $L^2$ and thus \as{}
\cite[Lemma 4.16]{Kallenberg}.

We can also study the joint distribution for several $m$. In particular,
\eqref{terminaliam} holds jointly for all $m\le n$, and thus 
\eqref{mormm} and \eqref{wsumm} hold jointly for all $m$.

Moreover, 
the conditional distribution of $M_n\mmx{m+1}$ 
given $M_n\mmx1,\dots,M_n\mm$
equals the
distribution of the maximum of $n-m$ \iid{} Exp(1) \rv{s}, conditioned on
this maximum being at most $M_n\mm$. In particular,
$M_n\mmx1,M_n\mmx2,\dots,M_n\mmx n$ form a Markov chain.
It is easy to see that this holds also in the limit as \ntoo. Thus
$W\mmx1,W\mmx2,\dots$ is a Markov chain, and the conditional distribution of
$W\mmx{m+1}$ given $W\mmx1,\dots,W\mm$ equals the distribution of 
$(W\mid W\le W\mm)$, where $W$ is a Gumbel variabel independent of $W\mm$.
Explicitly, for $x\le y$,
\begin{equation}
  \P\bigpar{W\mmx{m+1}\le x\mid W\mm=y}
=\P(W\le x\mid W\le y) = \exp\bigpar{e^{-y}-e^{-x}}.
\end{equation}
Note that this does not depend on $m$, so the Markov chain is
homogeneous. (The chain $M_n\mmx1,\dots,M_n\mmx n$ is not.)

This Markov chain was used by \citet{Fristedt} to describe the asymptotic
distribution of sizes the largest parts in a random partition (after
suitable normalization); the largest parts have the same asymptotic distribution
as the largest in a sequence of \iid{} exponential \rv{s}. 

The Markov chain becomes simpler if we transform to $e^{-W\mmx{m}}$.
The conditional distribution of $e^{-M_n\mmx{m+1}}$ given 
$e^{-M_n\mmx{1}},\dots,e^{-M_n\mmx{m}}$ equals the distribution of the
minimum $e^{-M'_{n-m}}$ 
of $n-m$ independent uniform random variables conditioned on this minimum being
at least $e^{-M_n\mmx{m}}$; in the limit it follows that
the conditional distribution of $e^{-W\mmx{m+1}}$ given 
$e^{-W\mmx{1}},\dots,e^{-W\mmx{m}}$ equals the distribution of 
$\bigpar{e^{-W'}\mid e^{-W'}\ge e^{-W\mm}}$, where $W'$ is a copy of $W$
independent of $W\mm$. Since $e^{-W'}\eqd e^{-W} \sim \Exp(1)$, it follows
that, conditionally given $W\mmx1,\dots,W\mm$,
\begin{equation}
e^{-W\mmx{m+1}}\eqd\bigpar{T\mid T\ge e^{-W\mm}}\eqd T+e^{-W\mm}  ,
\end{equation}
where $T\sim\Exp(1)$ is independent of $W\mm$. Consequently, the sequence
$e^{-W\mmx{1}},e^{-W\mmx{2}},\dots$ has the same distribution
as the sequence of partial sums of the \iid{} $\Exp(1)$ sequence
$T_1,T_2,\dots$:
\begin{equation}
\bigpar{ e^{-W\mmx{1}},e^{-W\mmx{2}},\dots}
\eqd
\bigpar{T_1,T_1+T_2,\dots}.
\end{equation}

In particular, this shows again that 
$e^{-W\mm}\sim\gG(m)$ and thus
$W\mm\eqd -\log \gG_m$.
 
The same asymptotic distributions $W\mm$ appear for the 
largest variables in many
other situations, see \eg{} \cite[Sections 2.2--2.3]{LLR}.
\end{exx}

\begin{exx}{Logistic distribution}\label{Elogi}
Let $\tW$ have the \emph{logistic distribution} with distribution function
$e^x/(e^x+1)$, or equivalently
\begin{equation}
  \P(\tW>x)=\frac1{e^x+1},
\qquad -\infty< x<\infty.
\end{equation}
By differentiation, the density function is
\begin{equation}\label{logidf}
  \frac{e^x}{(e^x+1)^2}
=
  \frac{1}{(e^{x/2}+e^{-x/2})^2}
=\frac1{4\cosh^2(x/2)}.
\end{equation}
If $\tP_1$ has the shifted Pareto distribution with density $(x+1)\qww$,
$x>0$, see Example 3.14, then 
$\P(\tP_1>e^x)=(e^x+1)\qw=\P(\tW>x)$
and thus
\begin{equation}
  \tW\eqd\log\tP_1.
\end{equation}
As a consequence, by Example 3.14,
$\tW$ has \mgf{} \ogt{} with
\begin{equation}
  \E e^{s\tW} =\E \tP_1^s
=\gG(1-s)\gG(1+s)
=\frac{\pi s}{\sin\pi s},
\qquad -1<\Re s<1.
\end{equation}
Equivalently, $\tW$ has the \chf
\begin{equation}\label{logichf}
  \E e^{\ii t\tW}
=\gG(1-\ii t)\gG(1+\ii t)
=\frac{\pi t}{\sinh \pi t}.
\end{equation}
We have $\rho_+=1$, $\rho_-=-1$, $\gam=2$, $\gam'=0$, $\gd=1$, $\gk=0$,
$C_1=2\pi$.

One way the logistic distribution appears is as the
symmetrization of the Gumbel distribution.
Let $W$ and $W'$ be \iid{} with the Gumbel
distribution (3.26), see Examples 3.19, 
and consider $W-W'$, which by (3.35) 
has the \mgf, 
for $ -1<\Re s<1$,
\begin{equation}
  \E e^{s(W-W')} =\E e^{sW}\E e^{-sW} 
=\gG(1-s)\gG(1+s)
=\E e^{s\tW}.
\end{equation}
Hence,
\begin{equation}\label{twww}
  \tW\eqd W-W'.
\end{equation}

By \eqref{twww} and \eqref{wsum} we further have the representation
\begin{equation}\label{twsum}
\tW\eqd
  \sum_{j\neq0} \frac{1}{j}(T_j-1)
=  \sum_{j=1}^\infty \frac{1}{j}(T_j-T_{-j})
\end{equation}
where $T_j$, $j\in\bbZ$, are \iid{} with the distribution $\Exp(1)$.
Since $T_j-T_{-j}$ has the \mgf{}
\begin{equation}
  \E e^{s(T_j-T_{-j})}
=\E e^{sT_1}\E e^{-sT_1}
=\frac1{1-s}\frac1{1+s}
=\frac1{1-s^2},
\qquad  -1<\Re s<1,
\end{equation}
\eqref{twsum} is equivalent to
\begin{equation}
  \frac{\pi s}{\sin \pi s} = \prod_{j=1}^\infty \frac{1}{1-s^2/j^2},
\end{equation}
which is a version of the product formula for $\sin$
\cite[4.3.89]{AS}
\begin{equation}
{\sin z} = z\prod_{j=1}^\infty \Bigpar{1-\frac{z^2}{j^2\pi^2}}.
\end{equation}

The random variable with the \chf{} $t/\sinh t$, and thus the distribution of
$\tW/\pi$, is studied by \citet{PitmanYor} (there denoted $\hat S_1$); 
among other things, they give the following construction:
Let $B(t)$ be a standard Brownian motion and let $T$ be the stopping time
when an independent standard 3-dimensional Brownian motion 
hits the
unit sphere in $\bbR^3$. Then
\begin{equation}
B(T)\eqd \tW/\pi.  
\end{equation}

The random variable $\tW/\pi$ (or $\tW$, depending on the choice of
normalization) appears also as the asymptotic distribution of the rank of a
random partition, see
\cite{SJpartitions}.
\end{exx}

\begin{exx}{Discriminants and Selberg's integral formula}
For a vector $(x_1,\dots,x_n)$ of real (or complex) numbers, define
\begin{equation}
  \gD\xxn\=\prod_{1\le i<j\le n} (x_j-x_i).
\end{equation}
Thus 
$\gD\xxn^2$ is the discriminant of the monic polynomial with roots
$x_1,\dots,x_n$.
Furthermore, $\gD\xxn$ is the well-known value of
the Vandermonde  determinant $\det\bigpar{x_i^{j-1}}_{i,j=1}^n$
(which apparently was never considered by Vandermonde, see \cite{MacTutor}).

Selberg \cite{Selberg} proved the following integral formula, for
$n\ge2$ and $\Re\ga>0$, $\Re\gb>0$, 
$\Re s>\max\set{-1/n, -\Re\ga/(n-1), -\Re\gb/(n-1)}$,
\begin{multline}\label{selberg}
  \intoi\dotsi\intoi |\gD\xxn|^{2s}\prodin x_i^{\ga-1}(1-x_i)^{\gb-1}\dd
  x_1\dotsm \dd x_n
\\
=\prod_{j=1}^{n}\frac{\gG\bigpar{\ga+(j-1)s}\gG\bigpar{\gb+(j-1)s)\gG(1+js}}
{\gG\bigpar{\ga+\gb+(n+j-2)s}\gG(1+s)}.
\end{multline}
(For applications of this formula, see \eg{} \cite{ForresterW} and
\cite{AndersonGZ}.) 
This leads to the following probabilistic interpretaions, see
\citet{LuRichards}.

For real $\ga,\gb>0$, 
let $\XXnx$ be \iid{} random variables with the Beta distribution
$\B(\ga,\gb)$; then
\eqref{selberg} can equivalently be written as the
expectation 
\begin{multline}\label{selberg2}
\E |\gD\XXn|^{2s}
\\
=\prod_{j=1}^{n}\frac
 {\gG(\ga+\gb)\gG\bigpar{\ga+(j-1)s}\gG\bigpar{\gb+(j-1)s}\gG(1+js)}
 {\gG(\ga)\gG(\gb)\gG\bigpar{\ga+\gb+(n+j-2)s}\gG(1+s)}.
\end{multline}
for $\Re s>\max\set{-1/n, -\Re\ga/(n-1), -\Re\gb/(n-1)}$.
This shows that 
$\gD\XXn^2$ has 
moments of Gamma type.
We have $\gam=\gam'=0$, $\gd=1-\ga-\gb-n/2$,
$\rho_+=\infty$ and $\rho_-=\max\set{-1/n, -\ga/(n-1), -\gb/(n-1)}$
(for $n\ge2$).

Equivalently, \eqref{selberg2} shows that the absolute value  $|\gD\XXn|$ has 
moments of Gamma type. 
In this case, see Remark 2.8,
$\gam=\gam'=0$, $\gd=1-\ga-\gb-n/2$,
$\rho_+=\infty$ and $\rho_-=2\max\set{-1/n, -\ga/(n-1), -\gb/(n-1)}$
(for $n\ge2$).

Note that for $n=1$, $\gD(X_1)=1$ is trivial, so the simplest non-trivial
case is $n=2$, when \eqref{selberg2}
says that if $X_1,X_2\sim\B(\ga,\gb)$ are independent, then
\begin{equation}
\E |X_1-X_2|^{2s}
=
\frac{\gG(\ga+\gb)^2}{\gG(\ga)\gG(\gb)}
\cdot
\frac
 {\gG\xpar{\ga+s}\gG\xpar{\gb+s}\gG(1+2s)}
 {\gG\xpar{\ga+\gb+s}\gG\xpar{\ga+\gb+2s}\gG(1+s)}.
\end{equation}

We obtain further results by taking suitable limits above,
\cf{} \cite{AndersonGZ}.
First, note that if $X_j\sim\B(\ga,\gb)$, then 
$\gb X_j\dto Y_j\sim\gG(\ga)$ as $\gb\to\infty$.
(For example by the method of moments, see (3.6) and (3.1).)
By taking limits in \eqref{selberg2}, using the facts that
\begin{equation}\label{ib1}
\gD(\gb X_1,\dots,\gb X_n)=\gb^{n(n-1)/2}\gD\XXn  
\end{equation}
and
\begin{equation}\label{ib2}
\gb^{-a}\frac{\gG(\gb+a)}{\gG(\gb)}\to1 
\qquad \text{as $\gb\to\infty$, for every fixed $a$},  
\end{equation}
it follows that if $\YYnx\sim\gG(\ga)$ are \iid{},
then
\begin{equation}\label{selbergY}
\E |\gD\YYn|^{2s}
=\prod_{j=2}^{n}\frac
 {\gG\bigpar{\ga+(j-1)s}\gG(1+js)}
 {\gG(\ga)\gG(1+s)}
\end{equation}
for $ \Re s>\max\set{-1/n,-\ga/(n-1)}$. 
Thus $\gD\YYn^2$ has moments of Gamma type, with 
$\gam=\gam'=n^2-n$, $\gd=(n-1)(\ga-1/2)$, $\rho_+=\infty$ and
$\rho_-=\max\set{-1/n,-\ga/(n-1)}$. 
In particular, $n=2$ yields
\begin{equation}
\E |Y_1-Y_2|^{2s}
=\frac
 {\gG\bigpar{\ga+s}\gG(1+2s)}
 {\gG(\ga)\gG(1+s)}.
\end{equation}
(For $\ga=1$, when $Y_1,Y_2\sim\Exp(1)$, this is an immediate consequence of
the fact that $|Y_1-Y_2|\sim\Exp(1)$ by the lack of memory in the
exponential distribution.)

Secondly, taking $\gb=\ga$, if $X_j\sim\B(\ga,\ga)$, then
$\sqrt{8\ga}(X_j-1/2)\dto Z_j\sim \N(0,1)$
as $\ga\to\infty$.
By taking limits in \eqref{selberg2}, using \eqref{ib1}--\eqref{ib2}
and the translation invariance
\begin{equation}\label{ib3}
\gD(X_1+a,\dots,X_n+a)=\gD\XXn  ,
\end{equation}
it follows that if $\ZZnx\sim\N(0,1)$ are \iid{},
then
\begin{equation}\label{selbergZ}
\E |\gD\ZZn|^{2s}
=\prod_{j=2}^{n}\frac
 {\gG(1+js)}
 {\gG(1+s)},
\qquad 
 \Re s>{-1/n}. 
\end{equation}
(This also follows by letting $\ga\to\infty$ in \eqref{selbergY}, 
using $(Y_\ga-\ga)/\sqrt\ga\dto\N(0,1)$, which for integer $\ga$ is just the
central limit theorem for $\gG(1)=\Exp(1)$.)
Thus $\gD\ZZn^2$ has moments of Gamma type, with 
$\gam=\gam'=n(n-1)/2$, $\gd=0$, $\rho_+=\infty$ and
$\rho_-=-1$.
Using the multiplication formula (A.5) for the Gamma function, 
\eqref{selbergZ} can be rewritten as
\begin{multline}\label{selbergZ2}
\E |\gD\ZZn|^{2s}
=
(2\pi)^{-n(n-1)/4}(n!)\qq
{\prod_{j=1}^n j^{js}}
\prod_{j=2}^{n}
\prod_{i=1}^{j-1}
{\gG(s+i/j)}
\\
=
\Bigpar{\prod_{j=1}^n j^{j}}^s
\prod_{1\le i<j\le n}
\frac {\gG(s+i/j)}{\gG(i/j)}
,
\qquad 
 \Re s>{-1/n}. 	
\end{multline}

The special case $n=2$ now just yields
\begin{equation}
\E |Z_1-Z_2|^{2s}
=\frac
 {\gG(1+2s)}
 {\gG(1+s)}
=\frac{2^{2s}}{\sqrt\pi}\gG(s+1/2),
\qquad 
 \Re s>-1/2,
\end{equation}
which is immediate because $Z_1-Z_2\sim\N(0,2)$, see (3.9).

The formulas for the moment imply some factorization formulas.
Thus, a  comparison between \eqref{selbergZ} and  Example 3.10 shows the
equality in distribution 
\begin{equation}
\gD\ZZn^2\eqd\prod_{j=2}^n S_{1/j}\qw, 
\end{equation}
where $S_{1/j}$ is stable with index
 $1/j$, and the variables are independent.
Similarly, \eqref{selbergZ2} and  (3.1) show the alternative factorization
 \cite{LuRichards}
\begin{equation}\label{selfact2}
  \gD\ZZn^2\eqd \prod_{j=1}^n j^j 
\prod_{1\le i<j\le n} G_{ij}
\end{equation}
with $G_{ij}\sim\gG(i/j)$ independent.

Similarly, 
\eqref{selbergY}, \eqref{selbergZ} and (3.1) yield
\begin{equation}
  \gD\YYn^2\eqd
\gD\ZZn^2\prod_{j=2}^n V_j^{j-1},
\end{equation}
where $V_j\sim\gG(\ga)$ are independent of each other and $\ZZnx$;
by \eqref{selfact2} this leads to a factorization of 
$\gD\YYn^2$ into independent Gamma variables.
Another such factorization is given by  \cite{LuRichards}, where also
similar factorizations of $\gD\XXn^2$ with $X\sim\B(\ga,\gb)$
are given for $n\le4$.
\end{exx}

\begin{exx}{Symmetric stable variables}\label{Esymmstab}
Consider a symmetric stable \rv{} $\bsa$ with \chf{}
$\gf(t)\=\E e^{\ii t \bsa}=e^{-|t|^\ga}$, where $0<\ga\le 2$.
(With this normalization, the L\'evy measure has density $c|x|^{-\ga-1}$,
where $c=\bigpar{-2\gG(-\ga)\cos\frac{\pi\ga}2}\qw$, see \eg{} \cite[Theorem
  3.3]{SJN12}.)

The (locally integrable function) $|x|^{s-1}$, where $0<\Re s<1$, has the
Fourier transform, in 
distribution sense, $c_s|x|^{-s}$ for a constant $c_s$
given by
\begin{equation}\label{cs}
  c_s=2^{s}\sqrt\pi\frac{\gG\parfrac{s}{2}}{ \gG\parfrac{1-s}{2}},
\end{equation}
see \eg{} \cite[Theorem IV.4.1]{Stein-Weiss}; 
this means that if
$\psi$ is in the 
Schwartz class $\cS$, then
\begin{equation}\label{fourier}
\intoooo |t|^{s-1}\hpsi(t)\dd t=c_s\intoooo|x|^{-s}\psi(x)\dd x,   
\end{equation}
where we
define the Fourier transform on $\cS$ by 
$\hpsi(x)\=\intoooo e^{\ii xt}\psi(t)\dd t$.
(The fact that the Fourier transform is of this type follows by a simple
homogeneity argument, and the value of $c_s$ then can be found by
considering the special case $\ga=2$ below.)

Since $\gf(t)\to0$ rapidly as $|t|\to\infty$,
$\bsa$ has a bounded and infinitely differentiable density function
$f(x)$; however, $f$ is not in $\cS$. Thus we regularize. Let $\eta(x)$ be
symmetric and 
infinitely differentiable with compact support and
$\heta(0)=\intoooo\eta(x)\dd x=1$, and define, for $\eps>0$,
$\eta_\eps(x)\=\eps\qw\eta(x/\eps)$, which has the Fourier transform
$\heta_\eps(x)=\heta(\eps x)$.
We then consider 
the product
$f_\eps(x)\=f(x)\heta_\eps(x)=f(x)\heta(\eps x)$, whose Fourier transform
is $\hf*\eta_\eps=\gf * \eta_\eps$; the function $f_\eps$ belongs to $\cS$,
and by applying 
\eqref{fourier} with $\psi=f_\eps$ and then letting $\eps\to0$, it follows
that
\eqref{fourier} holds with $\psi(x)=f(x)$ too, and thus, 
using \eqref{cs}
\begin{equation}
  \begin{split}
  \E|\bsa|^{-s}
&=\intoooo|x|^{-s}f(x)\dd x 
= c_s\qw  \intoooo|t|^{s-1}\gf(t)\dd t 
\\&
=2 c_s\qw  \intoo t^{s-1}e^{-t^\ga}\dd t 	
=2 c_s\qw \ga\qw \intoo u^{s/\ga-1}e^{-u}\dd u 	
\\&
=2 c_s\qw \ga\qw \gG(s/\ga)=
2^{1-s}\pi\qqw
\frac{\gG\parfrac{1-s}{2}\gG\parfrac{s}{\ga}}{\ga \gG\parfrac{s}{2}},
\\&
=
2^{-s}\pi\qqw
\frac{\gG\parfrac{1-s}{2}\gG\bigpar{1+\frac{s}{\ga}}}{\gG\bigpar{1+\frac{s}{2}}}.
  \end{split}
\end{equation}
We have proved this for $0<\Re s<1$, but by analytic continuation, it
extends to $-\ga<\Re s<1$, and thus
\begin{equation}\label{bsamom}
  \begin{split}
  \E|\bsa|^{s}
=
2^{s}
\frac{\gG\parfrac{1+s}{2}\gG\bigpar{1-\frac{s}{\ga}}}
 {\sqrt\pi\, \gG\bigpar{1-\frac{s}{2}}},
\qquad -1<\Re s<\ga.
  \end{split}
\end{equation}

Hence, $\bsa$ has \mogt. 
We have 
$\rho_+=\ga$ (except when $\ga=2$; then $\rho_+=\infty$) and $\rho_-=-1$;
furthermore, 
$\gam=1/\ga$, $\gam'=1-1/\ga$, $\gd=0$,
$\gk=\ga\qw\log\ga$, $C_1=\sqrt{4/\ga}$.

In the special case $\ga=2$, we have (with our choice of normalization)
$\bsa\eqd \sqrt2\, N$ with $N\sim\N(0,1)$, and thus \eqref{bsamom} is
equivalent to (3.9). (As said above, this yields a method to calculate $c_s$.)

In the special case $\ga=1$, $\bsx1$ has a Cauchy distribution with density
$1/(\pi(1+x^2))$, see \refE{Ecauchy}. 
In this case \eqref{bsamom} yields, using (A.3) and (A.6),
\begin{equation}\label{bsamom1}
  \begin{split}
  \E|\bsx1|^{s}
=
2^{s}
\frac{\gG\parfrac{1+s}{2}\gG\bigpar{1-s}}
 {\sqrt\pi\, \gG\bigpar{1-\frac{s}{2}}}
=
\frac{\gG\parfrac{1+s}{2}\gG\parfrac{1-s}2}
 {\pi}
=\frac1{\cos\frac{\pi s}{2}}
,
\quad -1<\Re s<1.
  \end{split}
\end{equation}
This is the same as \eqref{cauchy}, and also as
(3.13) with $n=1$. Indeed, it is well-known that
$\cT_1\eqd\bsx1$, \ie, $\cT_1$ has a Cauchy distribution.
\end{exx}

\begin{exx}{Products of Cauchy variables}\label{Ecauchy2} 
Let $X_1,X_2,\dots$ be \iid{} random variables with the Cauchy distribution
in \refE{Ecauchy}, and let $\Pi_k\=\prod_{1}^k X_i$ be the product of $k$ such
variables. 
(Note that $X_i\eqd X_i\qw$; thus \eg{} also $\Pi_2\eqd X_1/X_2$.)
It follows from  \refE{Ecauchy} that
$|\Pi_k|$ has moments of Gamma type
\begin{equation}
  \E |\Pi_k|^s=\frac1{\pi^k}
\gG\Bigpar{\frac 12+\frac s2}^k\gG\Bigpar{\frac 12-\frac s2}^k
=\frac1{\cos^k(\pi s/2)}
,
\qquad -1<\Re s<1.
\end{equation}
We have $\rho_+=1$, $\rho_-=-1$, $\gam=k$, $\gam'=0$, $\gd=0$,
$\gk=0$, $C_1=2^k$.

The density of $\Pi_1=X_1$ is $1/(\pi(1+x^2))$, and the density of
$\Pi_2=X_1X_2$ is
\begin{equation}
  \frac{2\log|x|}{\pi^2(x^2-1)}, \qquad -\infty<x<\infty,
\end{equation}
see \eg{} \citet{Pace}.
Formulas for the density of $\Pi_k$ for any integer $k\ge1$ are given by 
\citet{BFYor}.
\end{exx}

\begin{exx}{Generalized hyperbolic secant distribution}
The distribution of the L\'evy stochastic area $A$
in Example 3.20 is also known
as the \emph{hyperbolic secant distribution}, since both the density
function and the \chf{} are given by the hyperbolic secant $1/\cosh $ (up to
normalization constants).
This distribution is infinitely divisible, and thus, there exists a L\'evy
process $\hC_t$, $t\ge0$, such that $\hC_1=A$; consequently
$\hC_t$ has the \chf{}, \cf{} (3.36),
\begin{equation}\label{sech}
  \E e^{\ii s\hC_t} = \frac1{\cosh^t s}, \qquad s\in\bbR.
\end{equation}
The density is
\begin{equation}\label{sechdf}
  \frac{2^{t-2}}{\pi\gG(t)}\lrabs{\gG\parfrac{t+\ii x}{2}}^2,
\qquad x\in\bbR,
\end{equation}
see \eg{} \citet{PitmanYor} where many further results are given.

When $t=k$ is an integer, $\hC_k$ is the sum of $k$ independent copies of $A$,
so by Example 3.20, $\hC_k$ has \mgfogt, with 
\begin{equation}
  \E e^{s\hC_k} 
=\pi^{-k}\gG\Bigpar{\frac 12+\frac s\pi}^k\gG\Bigpar{\frac 12-\frac s\pi}^k
= \frac1{\cos^k s}, 
\qquad |\Re s|<\frac\pi2;
\end{equation}
we have $\rho_\pm=\pm\pi/2$, $\gam=2k/\pi$, $\gam'=0$,
$\gd=0$, $\gk=0$, $C_1=2^k$.
On the other hand, if $t$ is not an integer, then $\hC_t$ does not have \mgfogt,
since the \chf{} \eqref{sech} then cannot be extended to a meromorphic
function in $\bbC$.

The density \eqref{sechdf} is for $t=1$ 
\begin{equation}
  \frac{1}{2\cosh\frac{\pi x}2}
\end{equation}
as stated in Example 3.20, and for $t=2$
\begin{equation}\label{hc2}
  \frac{x}{2\sinh\frac{\pi x}2}.
\end{equation}
Similarly, for every integer $t=k\ge1$, the density \eqref{sechdf}
is a polynomial in $x$
divided by $\cosh(\pi x/2)$ ($k$ odd) or $\sinh(\pi x/2)$ ($k$ even); see
\citet{Harkness2} for explicit formulas. See also \cite{Holst} for an
application. 

Note that $\hC_k\eqd \frac2\pi\log|\Pi_k|$, where $\Pi_k$ is the product of
Cauchy variables in \refE{Ecauchy2}; this is an immediate consequence of the
case $k=1$ mentioned in \refE{Ecauchy}.

As a curiosity, we remark also that the distribution of $\hC_2$ is related
to the logistic distribution in \refE{Elogi} in the sense that the density
function of one distribution equals, up to constant factors and a rescaling,
the \chf{} of the other, see \eqref{logidf}, \eqref{logichf}, \eqref{sech},
\eqref{hc2}. 
In other words, the two density functions are
essentially the Fourier transforms of each other.
\end{exx}

\begin{exx}{Lamperti variables}\label{Elamperti}
Let $0<\ga<1$ and consider $L_\ga\=S_\ga/S'_\ga$ where $S_\ga,S'_\ga$ are two
independent copies of the 
  positive stable variable in Example 3.10; thus 
$\E e^{-tS_\ga}=\E e^{-tS'_\ga} =e^{-t^\ga}$, $t>0$.
By (3.16), $L_\ga$ has \mogt, using (A.6), 
\begin{equation}\label{lamperti}
  \begin{split}
	\E L_\ga^s 
&= \E S_\ga^s\E S_\ga^{-s}
= \frac{\gG(1-s/\ga)\gG(1+s/\ga)}{\gG(1-s)\gG(1+s)}
\\&
= \frac{\gG(s/\ga)\gG(1-s/\ga)}{\ga\gG(s)\gG(1-s)}
=\frac{\sin(\pi s)}{\ga\sin(\pi s/\ga)},
\qquad
-\ga < \Re s<\ga.
  \end{split}
\end{equation}
We have $\rho_\pm=\pm\ga$, 
$\gam=2\ga\qw-2$, $\gam'=\gd=\gk=0$, $C_1=1/\ga$, \cf{} Remarks
2.8 and 2.10. 

It is somewhat simpler to consider the power $L_\ga^\ga=M'_\ga/M_\ga$ where
$M_\ga,M'_\ga$ are \iid{} with the Mittag-Leffler distribution in Example
3.11.
By \eqref{lamperti}, \cf{} (3.17), 
\begin{equation}\label{lamperti2}
  \begin{split}
	\E (L_\ga^\ga)^s 
= \frac{\gG(1-s)\gG(1+s)}{\gG(1-\ga s)\gG(1+\ga s)}
=\frac{\sin(\pi \ga s)}{\ga\sin(\pi s)},
\qquad -1<\Re s<1.
  \end{split}
\end{equation}
We now have $\rho_\pm=\pm1$, 
$\gam=2-2\ga$, $\gam'=\gd=\gk=0$, $C_1=1/\ga$, \cf{} Remark 2.9.

The density of $L_\ga^\ga$ can be found by Fourier inversion, see \eg{}
\cite[p.~445]{Zolotarev:Mellin}, and can be written as
\begin{equation}\label{lamperti2df}
  \frac{\sin(\pi\ga)}{\pi\ga}\frac{1}{x^2+2\cos(\pi\ga)x+1},
\qquad x>0.
\end{equation}
Consequently, the density of $L_\ga$ is
\begin{equation}
  \frac{\sin(\pi\ga)}{\pi}\frac{x^{\ga-1}}{x^{2\ga}+2\cos(\pi\ga)x^\ga+1},
\qquad x>0.
\end{equation}

The random variable $L_\ga$ was studied (at least implicitly) by
\citet{Lamperti}, and is therefore called a \emph{Lamperti} variable by
\citet{James}, where also further references are given.

In the special case $\ga=1/2$, \eqref{lamperti2} simplifies to $1/\cos(\pi
s/2)$, so $L_{1/2}\qq$ is the absolute value of a Cauchy variable, see
\refE{Ecauchy}, which also follows directly from \eqref{lamperti2df}.

\citet{KotzO} defined, for $0<\ga<\gb\le2$,  a random variable
$Y_{\ga,\gb}=\bigpar{L_{\ga/\gb}^{\ga/\gb}}^{1/\ga}=L_{\ga/\gb}^{1/\gb}$.
(The defined $Y_{\ga,\gb}$ by giving its density function; that the definitions
are equivalent follows from \eqref{lamperti2df}.) By \eqref{lamperti} or
\eqref{lamperti2}, $Y_{\ga,\gb}$ has \mogt{}
\begin{equation}\label{lampertiKO}
  \begin{split}
\E Y_{\ga,\gb}^s
=\E L_{\ga/\gb}^{s/\gb} 
= \frac{\gG(1+s/\ga)\gG(1-s/\ga)}{\gG(1+s/\gb)\gG(1-s/\gb)}
=\frac{\gb\sin(\pi s/\gb)}{\ga\sin(\pi s/\ga)},
\qquad
-\ga < \Re s<\ga.
  \end{split}
\end{equation}
We have $\rho_\pm=\pm\ga$, $\gam=2\ga\qw-2$, $\gam'=\gd=\gk=0$, $C_1=1/\ga$.

\citet{James} also considers the more general
$X_{\ga,\gth}\=S_\ga/S_{\ga,\gth}$, for $\ga>0$ and $\gth>-\ga$,  where
$S_\ga$ and $S_{\ga,\gth}$ are 
independent, $S_\ga$ is a stable variable as above, and $S_{\ga,\gth}$ has a
distribution that is the same stable law tilted by $x^{-\gth}$, see
Remark 2.11. Thus $S_{\ga,\gth}$ has \mogt{} given by
\begin{equation}
\E S_{\ga,\gth}^s = \frac{\E S_\ga^{s-\gth}}{\E S_\ga^{-\gth}}
=
\frac{\gG(1+\gth)}{\gG(1+\gth/\ga)}
\frac{\gG(1-s/\ga+\gth/\ga)}{\gG(1-s+\gth)},
\qquad \Re s < \ga+\gth,
\end{equation}
and $X_{\ga,\gth}$ has \mogt{} given by,
for $ - \ga-\gth<\Re s < \ga$,
\begin{equation}
\E X_{\ga,\gth}^s = 
\E S_\ga^s \E S_{\ga,\gth}^{-s} 
=
\frac{\gG(1+\gth)}{\gG(1+\gth/\ga)}
\frac{\gG(1-s/\ga)\gG(1+s/\ga+\gth/\ga)}
{\gG(1-s)\gG(1+s+\gth)}.
\end{equation}
We have $\rho_+=\ga$, $\rho_-=-\ga-\gth$, 
$\gam=2(1/\ga-1)$, $\gam'=0$, $\gd=\gth(1/\ga-1)$, $\gk=0$.
\end{exx}

\begin{exx}{A generalized exponential distribution}\label{Etom}
  Let $\gb>0$ and let $V_\gb$ be a positive random variable with the density
  function 
  \begin{equation}
	\frac{1}{\gG\lrpar{1+\xfrac1\gb}}e^{-x^\gb},
\qquad x>0.
  \end{equation}
A simple change of variables verifies that this is a probability density
function, and more generally that,
for $\Re s>-1$,
\begin{equation}\label{tom}
  \E V_\gb^s = \frac1{\gG(1+1/\gb)}\int_0^\infty x^s e^{-x^\gb}\dd x
=\frac{\gG(s/\gb+1/\gb)}{\gb\gG(1+1/\gb)}
=\frac{\gG(s/\gb+1/\gb)}{\gG(1/\gb)}.
\end{equation}
$V_\gb$ thus has \mogt, with
$\rho_+=\infty$, $\rho_-=-1$, $\gam=\gam'=1/\gb$, $\gd=\frac1\gb-\frac12$,
$\gk=\frac1\gb\log\frac1\gb$.

Note that $\gb=1$ gives the exponential distribution in Example 3.2. In
general, the distribution of $V_\gb$ can be seen as a tilted version of the
Weibull distribution in Example 3.7.
\end{exx}

\begin{exx}{Linnik distribution}
The Linnik distribution \cite{Linnik} has \chf{}
\begin{equation}
  \frac1{1+|t|^\ga},
\end{equation}
where $0<\ga\le2$.
As shown by \citet{Devroye}, a random variable $X_\ga$
with this distribution is
easily constructed as
\begin{equation}
  X_\ga\=\bsa V_1^{1/\ga},
\end{equation}
where $\bsa$ is the symmetric stable  random variable in
\refE{Esymmstab}, $V_1$ has the exponential distribution $\Exp(1)$,
and these are independent.
  
By \eqref{bsamom} and (3.2),
for $-\min(\ga,1)<\Re s<\ga$,
\begin{equation}\label{linnik}
  \begin{split}
\E|X_\ga|^s=
  \E|\bsa|^{s} \E V_1^{s/\ga}
=
2^{s}
\frac{\gG\parfrac{1+s}{2}\gG\bigpar{1+\frac{s}{\ga}}\gG\bigpar{1-\frac{s}{\ga}}}
 {\sqrt\pi\, \gG\bigpar{1-\frac{s}{2}}}.
  \end{split}
\end{equation}
Hence $X_\ga$ has \mogt, with $\rho_+=\ga$, $\rho_-=-\min(\ga,1)$, 
$\gam=2/\ga$, $\gam'=1$, $\gd=1/2$, $\gk=0$, $C_1=\sqrt{2\pi}$.

More generally, \citet{Devroye} showed that if $0<\ga\le2$ and $\gb>0$, and
$\bsa$ 
is as above and $Y_\gb$ as in \refE{Etom} and independent of $\bsa$, then
\begin{equation}
  X_{\ga,\gb}\=\bsa V_\gb^{\gb/\ga}
\end{equation}
has \chf
\begin{equation}
  \frac1{(1+|t|^\ga)^{1/\gb}}.
\end{equation}
(This implies that $X_\ga=X_{\ga,1}$, and more generally every
$X_{\ga,\gb}$, 
is infinitely divisible, and that there is a \Levy{} process $\hat X_{\ga,t}$,
$t\ge0$, such that $\hat X_{\ga,t}\eqd X_{\ga,1/t}$ for all $t>0$.)

By \eqref{bsamom} and \eqref{tom},
for $-\ga/\gb<\Re s<\ga$,
\begin{equation}
  \begin{split}
\E|X_{\ga,\gb}|^s=
  \E|\bsa|^{s} \E V_\gb^{s\gb/\ga}
=
2^{s}
\frac{\gG\parfrac{s+1}{2}\gG\bigpar{1-\frac{s}{\ga}}
 \gG\bigpar{\frac{s}{\ga}+\frac1\gb}}
 {\sqrt\pi\, \gG\bigpar{\frac1\gb}\gG\bigpar{1-\frac{s}{2}} }.
  \end{split}
\end{equation}
Hence $X_{\ga,\gb}$ has \mogt, with $\rho_+=\ga$, $\rho_-=-\ga/\gb$, 
$\gam=2/\ga$, $\gam'=1$, $\gd=1/\gb-1/2$, $\gk=0$.

\citet{KotzO} showed that $X_\ga\eqd X_\gb Y_{\ga,\gb}$ where $Y_{\ga,\gb}$
is as in \refE{Elamperti} and independent of $X_\gb$; this follows also
directly from \eqref{lampertiKO} and \eqref{linnik}.
For the Linnik distribution see further \cite{KotzOH}.
\end{exx}

\newcommand\AAP{\emph{Adv. Appl. Probab.} }
\newcommand\JAP{\emph{J. Appl. Probab.} }
\newcommand\JAMS{\emph{J. \AMS} }
\newcommand\MAMS{\emph{Memoirs \AMS} }
\newcommand\PAMS{\emph{Proc. \AMS} }
\newcommand\TAMS{\emph{Trans. \AMS} }
\newcommand\AnnMS{\emph{Ann. Math. Statist.} }
\newcommand\AnnAP{\emph{Ann. Appl. Probab.} }
\newcommand\AnnPr{\emph{Ann. Probab.} }
\newcommand\CPC{\emph{Combin. Probab. Comput.} }
\newcommand\JMAA{\emph{J. Math. Anal. Appl.} }
\newcommand\RSA{\emph{Random Struct. Alg.} }
\newcommand\ZW{\emph{Z. Wahrsch. Verw. Gebiete} }
\newcommand\DMTCS{\jour{Discr. Math. Theor. Comput. Sci.} }

\newcommand\AMS{Amer. Math. Soc.}
\newcommand\Springer{Springer-Verlag}
\newcommand\Wiley{Wiley}

\newcommand\vol{\textbf}
\newcommand\jour{\emph}
\newcommand\book{\emph}
\newcommand\inbook{\emph}
\def\no#1#2,{\unskip#2, no. #1,} 
\newcommand\toappear{\unskip, to appear}

\newcommand\webcite{\url}
\newcommand\webcitesvante{\webcite{http://www.math.uu.se/~svante/papers/}}
\newcommand\arxiv[1]{\webcite{arXiv:#1.}}

\def\nobibitem#1\par{}

\end{document}